\documentclass[11pt,reqno]{amsart}
\usepackage{amsmath, amssymb, amsthm}
\usepackage{url}
\usepackage[breaklinks]{hyperref}
\usepackage{enumitem}
\allowdisplaybreaks
 \usepackage{autonum}

\usepackage{color}

\theoremstyle{plain}
\newtheorem{thm}{Theorem}[section]
\theoremstyle{plain}

\newtheorem{prop}[thm]{Proposition}
\newtheorem{cor}[thm]{Corollary}

\theoremstyle{definition}

\numberwithin{equation}{section}

\begin{document}

\title{Divisibility properties of weighted $k$-regular partitions}

\date{\today}

\author{Debika Banerjee and Ben Kane}
\address{Debika Banerjee\\ Department of Mathematics\\
Indraprastha Institute of Information Technology IIIT, Delhi\\
Okhla, Phase III, New Delhi-110020, India.} 
\email{debika@iiitd.ac.in}

\address{Ben Kane\\ Department of Mathematics\\
University of Hong Kong\\
Pokfulam, Hong Kong.} 
\email{bkane@hku.hk}

\thanks{2010 \textit{Mathematics Subject Classification.} Primary 05A17, 11P81, 11P83, Secondary 11F11.\\
\textit{Keywords and phrases.} Modular form, cusp form, $k$-Regular partitions, congruence.
}

 \begin{abstract}
We study a generalized class of weighted $k$-regular partitions defined by  
\[
\sum_{n=0}^{\infty} c_{k, r_1, r_2}(n) q^n  
= \prod_{n=1}^{\infty} \frac{(1 - q^{nk})^{r_1}}{(1 - q^n)^{r_2}},
\]  
which extends the classical $k$-regular partition function $b_k(n)$. We establish new infinite families of Ramanujan-type congruences, divisibility results, and positive-density prime sets for which $c_{k, r_1, r_2}(n)$ vanishes modulo a given prime. 
\end{abstract}

\maketitle

\section{Introduction}

Partition theory plays a central role in number theory and combinatorics, with origins in the classical work of Euler, Ramanujan, and Hardy. Partitions of $n$ are non-increasing sequences of positive integers (parts) which sum to $n$. Among the various classes of partitions, the family of \emph{$k$-regular partitions}---those in which no part is divisible by a fixed integer $k \geq 1$---has been extensively studied. The generating function for the number of such partitions of an integer $n$, denoted by $b_k(n)$, is given by
\begin{align}\label{k11}
\sum_{n=0}^{\infty} b_k(n)\, q^n = \prod_{n=1}^{\infty} \frac{1 - q^{nk}}{1 - q^n}.
\end{align}
This family generalizes the classical partition function $p(n)$, which counts the number of partitions of $n$, since $b_k(n) = p(n)$ for $k > n$. and arises naturally in the theory of modular forms, $q$-series, and congruences in arithmetic functions. Moreover, when $k=p$ is prime, $b_p(n)$ coincides with the number of irreducible $p$-modular representations of the symmetric group $S_n$~\cite{JK}.

Pioneering results on the congruence properties of $k$-regular partitions were established by Gordon and Hughes~\cite{GordonHughes1}, and further extended in later works by Gordon and Ono~\cite{GordonHughes2}, Hirschhorn and Sellers~\cite{H}, and Lovejoy~\cite{JL}. These investigations, together with Martin's theory of multiplicative $\eta$-quotients~\cite{Mar}, have revealed deep connections between partition functions and modular forms. The modular approach has been developed extensively in the works of Murty~\cite{Mur}, and Ono~\cite{Ono2004}, with modern surveys such as Ono and Webb~\cite{OnoWebb} offering a comprehensive perspective. In particular, Mahlburg~\cite{Mahlburg2005} extended Ramanujan-type congruences using $p$-adic modular forms, opening the way for further density and divisibility results.

The arithmetic and asymptotic behavior of $k$-regular partitions have also been studied from analytic and Galois-theoretic perspectives. Serre~\cite{Serre1973} showed that many partition functions modulo primes are modular and exhibit congruences on positive density sets. Ahlgren and Ono~\cite{AhlgrenOno2001} investigated congruences and density results for partition functions using cusp forms and Galois representations. More recently, Zheng~\cite{Zheng2022} examined the divisibility and distribution of the $3$-regular and $5$-regular partition functions $b_3(n)$ and $b_5(n)$, investigating both regular and irregular behaviors in its residue classes modulo primes.

Beyond the standard $k$-regular partitions, one can consider \emph{generalized or weighted partition functions} that extend the form of~\eqref{k11}. Let $r_1 \geq 1$ and $r_2 \geq 1$ be integers. Define the function $c_{k, r_1, r_2}(n)$ via the generating function
\begin{align}\label{defico}
\sum_{n=0}^{\infty} c_{k, r_1, r_2}(n) q^n 
= \prod_{n=1}^{\infty} \frac{(1 - q^{nk})^{r_1}}{(1 - q^n)^{r_2}}.
\end{align}
If $r_2\geq r_1$, then $c_{k, r_1, r_2}(n)$ counts the number of weighted partitions of $n$ in which parts divisible by $k$ appear in $r_2 - r_1$ colors, and all other parts appear in $r_2$ colors. Clearly, for $r_1 = r_2 = 1$, we recover $b_k(n)$.

These generalized partition functions arise in the study of combinatorial generating functions, identities of the Rogers--Ramanujan type, and in asymptotic enumeration problems. Although we assume that $r_2\geq r_1$ throughout this paper, we note that  $c_{k,k,1}(n)$ counts the number of $k$-core partitions of $n$, which were shown by Ono \cite{Onotcore1,Onotcore2} and Granville--Ono \cite{GranvilleOno} to be positive for $k\geq 4$ (i.e., $c_{k,k,1}(n)>0$ for $k\geq 4$ and $n\in\mathbb{N}$). Taking $r_1=r_2=r$, one recovers the $k$-th infinite Borwein product $\prod_{n\geq 1} \frac{1-q^n}{1-q^{kn}}$ taken to the power $-r$; if one allows negative $r$, then signs of these Fourier coefficients were considered as special cases of results of Schlosser and Zhou \cite{SchlosserZhou}.  

In this paper, we extend the study of $5$-regular partitions by analyzing the arithmetic and asymptotic behavior of the generalized partition function $c_{k, r_1, r_2}(n)$. Our aim to is
build upon the results of Zheng ~\cite{Zheng2022} in order to obtain congruences
satisfied by the weighted $k$-regular partition function $c_{k, r_1, r_2}(n)$ modulo primes.

\begin{thm}\label{thm:main}
Let \( p \) be a prime and \( M \geq 1 \) an odd integer satisfying \( p \geq M \).  
Let \( m \geq 5 \) be a sufficiently large prime, such that $\frac{(p-1)(m-1)}{4}$ is even and let \( r \geq 1 \) be an integer.  
Let \( c_{p, r, Mr}(n) \) be as defined in~\eqref{defico}, and set \( s = \gcd(r, 24) \) and \( d = \gcd\!\left( \frac{24}{s},\, p - M \right) \).  
Then there exists a set of primes \( \ell \) of positive density such that
\[
c_{p, r, Mr}\!\left( \frac{d m n \ell - r(p - M)}{24} \right) \equiv 0 \pmod{m},
\]
for all integers \( n \) with \( \gcd(n, \ell) = 1 \).
\end{thm}
In particular, if we put \( r = M = 1 \), we obtain the following corollary.

\begin{cor}
Let \( p \) be a prime and let \( m \geq 5 \) be a sufficiently large prime.  
Let \( b_{k}(n) \) be as defined in~\eqref{k11}, and set \( d = \gcd(24,\, p - 1) \).  
Then there exists a set of primes \( \ell \) of positive density such that
\[
b_{p}\!\left( \frac{d m n \ell - (p - 1)}{24} \right) \equiv 0 \pmod{m},
\]
for all integers \( n \) with \( \gcd(n, \ell) = 1 \).
\end{cor}
In particular, when $p\in\{3,5\}$, we recover \cite[Theorems 2 and 3]{Zheng2022}. The proof mostly follows the argument of Zheng \cite{Zheng2022}, with one major technical difficulty. In Zheng's case, the minimal order of vanishing coming from \cite[Theorem 5]{Zheng2022} immediately implies that a certain ratio is a cusp form because of the level. In our more general setting, this does not immediately follow because the level goes to infinity. As a result, we need to more carefully analyze the order of vanishing at all cusps of the image of a cusp form under the Hecke operator, showing that the main term vanishes (see the calculations leading up to \eqref{eqn:gcuspidal}). The other steps follow from Zheng's argument in \cite{Zheng2022} and known results in the theory of modular forms.

The paper is organized as follows. In Section \ref{sec:prelim}, preliminaries about eta-products and modular forms are recalled. In Section \ref{sec:modularconstruction}, we construct certain cusp forms which play an important role in our proofs. Finally, we compute the order of vanishing at all cusps under the Hecke operator and prove Theorem \ref{thm:main} in Section \ref{sec:main}.

\section*{Acknowledgements}

The research of the first author was funded by the SERB under File No. MTR/2023/000837 and CRG/2023/002698. 
The research of the second author was supported by grants from the Research Grants Council of the Hong Kong SAR, China (project numbers HKU 17314122, HKU 17305923).


\section{Preliminaries
}\label{sec:prelim}

In this section, we collect foundational results and notation related to eta-quotients and $k$-regular partitions that will be used throughout the paper.

\subsection{Modular forms}
Recall that, for $\kappa,N\in\mathbb{N}$ and a character $\chi$ modulo $N$, a holomorphic function $f$ from the complex upper half-plane $\mathbb{H}$ to $\mathbb{C}$ is called a \begin{it}modular form\end{it} of weight $\kappa$ and Nebentypus character $\chi$ if for every $\gamma=\left[\begin{smallmatrix}a&b\\c&d\end{smallmatrix}\right]\in\Gamma_0(N)$ it satisfies 
\begin{equation}\label{eqn:modulardef}
\chi(d) f(\tau)=f|_{\kappa}\gamma(z):=(c\tau+d)^{-\kappa} f\left(\frac{a\tau+b}{c\tau+d}\right),
\end{equation}
and $f|_\kappa\gamma(\tau)$ grows at most polynomially as $\tau\to i\infty$ for every $\gamma\in\operatorname{SL}_2(\mathbb{Z})$. More generally, for $\kappa\in\mathbb{Z}$ we say that a function $f:\mathbb{H}\to\mathbb{C}$ satisfies weight $k$ \begin{it}modularity\end{it} on $\Gamma_0(N)$ with character $\chi$ if \eqref{eqn:modulardef} holds. We call the equivalence classes in $\Gamma\backslash (\mathbb{Q}\cup\{i\infty\})$ the \begin{it}cusps\end{it} of $\Gamma$, and sometimes abuse notation to call elements of $\mathbb{Q}\cup\{i\infty\}$ the cusps of $\Gamma$. The condition that $f|_\kappa\gamma(z)$ grows at most polynomially as $\tau\to i\infty$ may then be considered a growth condition of $f$ towards the cusp $\gamma(i\infty)=\frac{a}{c}$. For a cusp $\alpha=\gamma(i\infty)$ and $q:=e^{2\pi i \tau}$, we call 
\[
f|_\kappa\gamma(\tau)=\sum_{n\gg -\infty} a_{f,\alpha}(n) q^{\frac{n}{M}}
\]
the \begin{it}Fourier expansion\end{it} of $f$ at $\alpha$. Here $M$ is called the \begin{it}cusp width\end{it} of $\alpha$ and may be chosen minimally so that $T^M\in \gamma^{-1}\Gamma \gamma$. We omit $\alpha$ in the notation when it is $i\infty$ and sometimes omit $f$ when it is clear from context.

\subsection{Eta-Quotients and Modularity}

Let $\eta(z)$ denote the \textit{Dedekind eta-function}, defined by
\[
\eta(\tau) := q^{1/24} \prod_{n=1}^\infty (1 - q^n), \quad \text{where } q = e^{2\pi i \tau}, \ \Im(z) > 0.
\]
An \textit{eta-quotient} is a function of the form
\[
f(\tau) = \prod_{\delta \mid N} \eta^{r_\delta}(\delta \tau),
\]
for integers $r_\delta\in\mathbb{Z}$ indexed by the positive divisors $\delta$ of some fixed $N \in \mathbb{N}$. The modular properties of such functions are described by the following result.

\begin{thm}[Gordon–Hughes–Newman]\label{thm1}
Let \( f(\tau) = \prod_{\delta \mid N} \eta^{r_\delta}(\delta \tau) \) be an eta-quotient satisfying the conditions
\begin{align}\label{one}
\sum_{\delta \mid N} \delta r_\delta \equiv 0 \pmod{24}, \quad \text{and} \quad
\sum_{\delta \mid N} \frac{N}{\delta} r_\delta \equiv 0 \pmod{24}.
\end{align}
Then \( f(\tau) \) transforms under $\Gamma_0(N)$ as
\begin{align}\label{three}
f(A\tau) = \chi(d) (c\tau + d)^{\kappa} f(\tau),
\end{align}
for all \( A = \begin{bmatrix} a & b \\ c & d \end{bmatrix} \in \Gamma_0(N) \), where the weight is
\[
\kappa = \frac{1}{2} \sum_{\delta \mid N} r_\delta,
\]
and the character $\chi$ is given by
\[
\chi(d) = \left( \frac{(-1)^{\kappa} s}{d} \right), \quad \text{with} \quad s = \prod_{\delta \mid N} \delta^{r_\delta}.
\]
In other words, $f$ satisfies weight $\kappa$ modularity on $\Gamma_0(N)$ with character $\chi$.
\end{thm}
The behavior of eta-quotients at the cusps can be described explicitly using the following result.

\begin{thm}[Ligozat]\label{thm2}
Let \( f(\tau) = \prod_{\delta \mid N} \eta^{r_\delta}(\delta \tau) \) be an eta-quotient satisfying the modularity conditions \eqref{one}. Let $d \mid N$, and let $c$ be a positive integer with $\gcd(c,d) = 1$. Then the order of vanishing of $f(\tau)$ at the cusp $c/d$ is
\[
\frac{1}{24} \cdot \frac{N}{d \cdot \gcd(d, N/d)} \sum_{\delta \mid N} \frac{r_\delta}{\delta} \cdot \gcd(d, \delta)^2.
\]
\end{thm}
We need a famous theorem of Serre about congruences for the Hecke operators.
\begin{thm}[J.-P. Serre]\label{thm:serre}
The set of primes \( \ell \equiv -1 \pmod{Nm} \) such that
\[
f \mid T(\ell) \equiv 0 \pmod{m},
\]
for each \( f(\tau) \in S_k(\Gamma_0(N), \psi)_m \) has positive density, where \( T(\ell) \) denotes the usual Hecke operator acting on \( S_k(\Gamma_0(N), \psi) \).
\end{thm}

Before proving the key result, we record some useful facts about modular forms and their Fourier coefficients.  
For a more detailed overview, the reader is referred to~\cite{NK}.

\begin{prop}\label{prop1}
Let \( f(\tau) = \sum_{n=0}^\infty a(n) q^n \) be a modular form in \( M_k(\Gamma_0(N), \psi) \).  
\begin{enumerate}
    \item For any positive integer \(t\), the function
    \[
    f(t\tau) = \sum_{n=0}^\infty a(n)\, q^{tn},
    \]
is the Fourier expansion of a modular form in \( M_k(\Gamma_0(tN), \psi) \).

    \item For any prime \(p\), the function
    \[
    f(\tau)\,\big|\,T(p) := \sum_{n=0}^\infty \Big( a(pn) + \psi(p)\,p^{k-1} a\!\left(\tfrac{n}{p}\right) \Big) q^n,
    \]
    is the Fourier expansion of a modular form in \( M_k(\Gamma_0(N), \psi) \).
\end{enumerate}
Moreover, both assertions remain valid when \( M_k(\Gamma_0(N), \psi) \) is replaced by \( S_k(\Gamma_0(N), \psi) \).

Here \(T(p)\) denotes the Hecke operator associated with the prime \(p\).  
In particular, one has the congruence
\[
f(\tau)\,\big|\,T(p) \equiv f(\tau)\,\big|\,U(p) \pmod{p},
\]
in \( M_k(\Gamma_0(N), \psi)_{\mathbb{F}_p} \), where the operator \(U(p)\) acts on a Fourier series by
\[
\left(\sum_{n=0}^\infty a(n)\,q^n\right)\big|\,U(p) := \sum_{n \equiv 0 \,(\mathrm{mod}\,p)} a(n)\, q^{n/p}.
\]
This criterion enables one to verify congruences between modular forms by means of a finite computation.
\end{prop}

\subsection{Functional Equation and Modular Transformations}

The Dedekind eta function satisfies the transformation law
\begin{align}\label{fnal equation}
\eta\left(-\frac{1}{\tau}\right) = (-i\tau)^{1/2} \eta(\tau),
\end{align}
which plays a central role in the modular transformation behavior of eta-products.

To study modular transformations more precisely, especially under $\mathrm{SL}_2(\mathbb{Z})$, we often work with exponential parametrizations of the upper half-plane. Let $w \in \mathbb{C}$ with $\Re(w) > 0$, and define
\[
q = \exp\left( \frac{2\pi i}{k}(h + iw) \right),
\]
where $h, k \in \mathbb{N}$, $0 \leq h < k$, and $\gcd(h, k) = 1$. Let $h'$ satisfy $hh' \equiv -1 \pmod{k}$, and define
\[
q_1 = \exp\left( \frac{2\pi i}{k}\left(h' + iw^{-1}\right) \right).
\]

With this notation, the transformation of $\eta$ between the points $\frac{1}{k}(h + iw)$ and $\frac{1}{k}(h' + iw^{-1})$ is determined by a certain constant $w_{h,k}$ (known as the $\eta$-multiplier) as follows:
\begin{align}\label{whk}
\eta\left( \frac{1}{k}(h + iw) \right) 
&= w_{h,k}^{-1} e^{-\frac{\pi i}{4}} w^{-\frac{1}{2}} \, e^{\frac{\pi i}{12k}(h - h')} \, \eta\left( \frac{1}{k}(h' + iw^{-1}) \right), \\
\label{poa}
\eta\left( \frac{1}{k}(h' + iw^{-1}) \right) 
&= w_{h,k} \, e^{\frac{\pi i}{4}} w^{\frac{1}{2}} \, e^{-\frac{\pi i}{12k}(h - h')} \, \eta\left( \frac{1}{k}(h + iw) \right),
\end{align}
where
\begin{align}\label{poq}
w_{h,k} = \left( \frac{-h}{k} \right) 
\exp\left( -\pi i \left[ \frac{1}{4}(k - 1) + \frac{1}{12} \left( k - \frac{1}{k} \right)(h - h') \right] \right).
\end{align}
The above expressions are crucial in analyzing modular transformations of eta-quotients at non-trivial cusps.

\section{Modular Construction via Eta-Products Modulo \texorpdfstring{$m$}{m}}\label{sec:modularconstruction}

Let $m \geq 5$ be a prime. Also assume $(m, r_1)=(m, r_2)=1$. Define the eta-product
\begin{align}
    f_{p, m, r_1, r_2}(\tau) := \left(\frac{\eta(p\tau)^{r_1}}{\eta(\tau)^{r_2}}\right)
    \eta(pm\tau)^a \eta(m\tau)^b,
\end{align}
where $a, b$ are positive integers to be determined. Using the well-known congruence
\[
(1 - x^p)^m \equiv (1 - x^{pm}) \pmod{m},
\]
it follows that
\[
f_{p, m, r_1, r_2}(\tau) \equiv \eta(p\tau)^{am + r_1} \eta(\tau)^{bm - r_2} \pmod{m}.
\]

Our goal is to select $a$ and $b$ such that the resulting eta-quotient satisfies the conditions in Theorem~\ref{thm1}. These lead to the congruences
\begin{align}\label{eq1}
    m(pa + b) \equiv r_2 - pr_1 \pmod{24},
\end{align}
and
\begin{align}\label{eq2}
    m(a + pb) \equiv pr_2 - r_1 \pmod{24}.
\end{align}

\vspace{0.5cm}

\section{Congruences for the case \texorpdfstring{$r_1 = r$, $r_2 = M r$ with $M$ odd}{r1=r, r2=Mr with M odd}}\label{sec:main}

We note that, for $m$ sufficiently large (so that $bm\geq Mr$), if \eqref{eq1} and \eqref{eq2} are satisfied, then Theorem \ref{thm1} implies that 
\begin{align}\label{modularform}
\eta^{am + r}(p\tau)\, \eta^{bm - Mr}(\tau)
\in S_{\kappa}(\Gamma_0(p), \chi_p),
\end{align} with
\begin{align}\label{wt}
\kappa = \frac{(a + b)m + r(1 - M)}{2}.
\end{align}

Let $s = \gcd(r, 24)$, and write $r = sv$. We take $a = s a_1$, $b = s b_1$ for some integers $a_1$, $b_1$. Substituting into \eqref{eq1} and \eqref{eq2}, we obtain:
\begin{align}
    m(pa_1 + b_1) &\equiv v(M - p) \pmod{\frac{24}{s}}, \label{a1eq}\\
    m(a_1 + pb_1) &\equiv v(pM - 1) \pmod{\frac{24}{s}}. \label{b1eq}
\end{align}
Let $d = \gcd(p - M, \tfrac{24}{s})$. We now choose
\[
a_1 = p - m', \qquad b_1 = M m' - 1,
\]
so that
\[
m m' \equiv v \pmod{\frac{24}{ds}}.
\]
This choice implies that $m' \equiv m^{-1} v\equiv mv \pmod{\frac{24}{ds}}$.

We now examine the $q$-expansion of $f_{p, m,r,Mr}(\tau)$. From the definition, we have:
\begin{align}
    f_{p, m, r, Mr}(\tau) &= \sum_{n \geq 0} c_{p, r, Mr}(n)\, q^{\frac{24n + r(p - M)}{24}} \cdot q^{\frac{m(pa + b)}{24}} \cdot
    \prod_{n=1}^{\infty} (1 - q^{pmn})^a (1 - q^{mn})^b.
\end{align}

Applying the $U(m)$ operator and using the congruence $U(m) \equiv T(m) \pmod{m}$, we obtain
\begin{align}\label{MO}
    \sum_{n \geq 0} c_{p, r, Mr}(n)\, q^{\frac{24n + m(pa + b) + r(p - M)}{24}} \bigg| U(m)
    \equiv \frac{\eta^{am + r}(p\tau)\, \eta^{bm - Mr}(\tau)\big| T(m)}
    {\prod_{n=1}^{\infty}(1 - q^{pn})^a (1 - q^n)^b} \pmod{m},
\end{align}
where $T(m)$ denotes the usual Hecke operator acting on $S_k(\Gamma_0(p), \chi_p)$.

To isolate a specific subsequence, we extract coefficients for which
\[
m \mid \left( \frac{24}{d} n + \frac{r(p - M)}{d} \right),
\]
where $d = \gcd(p - M, 24)$. This yields:
\begin{align}
    \sum_{\substack{n \geq 0 \\ m \mid \left(\frac{24}{d}n + \frac{r(p - M)}{d} \right)}}
    c_{p, r, Mr}(n)\, q^{\frac{\frac{24}{d}n + \frac{r(p - M)}{d}}{(24/d)m}}
    &\equiv \frac{\eta^{am + r}(pz)\, \eta^{bm - Mr}(z)\big| T(m)}{\eta(pz)^a \eta(z)^b} \pmod{m}. \label{MM}
 \end{align}
 Let use denote the right hand side of \eqref{MM} by $g(q)$. This congruence provides the foundation for establishing congruence relations for the coefficients $c_{p, r, Mr}(n)$ modulo $m$ through the modular form $g(q)$.

\subsection{Modular Transformation Behavior of Eta Products}

Let $\Re(z) > 0$ throughout this section. We next determine the growth of the function $g(q)$ appearing on the right-hand side of \eqref{MM} towards all of the cusps. We first deal with the denominator of $g(q)$.

\noindent{\bf Case 1: If $p \nmid k$.}\\
Let $k \in \mathbb{N}$ with $\gcd(h_0, k) = 1$, and suppose $p \nmid k$. Set $h' = ph_0$ and $w = \frac{z}{p}$, so that $q = e^{\frac{2\pi i}{k} (h'+ i/w)}$. We choose $h \in \mathbb{Z}$ such that 
\[
h h' \equiv -1 \pmod{k}.
\]

Applying the transformation identity~\eqref{poa}, we obtain
\begin{align}
\eta\left( \frac{1}{k}\left(ph_0 + \frac{ip}{z} \right) \right)
&= \eta\left( \frac{1}{k}(h' + i w^{-1}) \right)\\
&= w_{h,k} e^{\frac{\pi i}{4}} \sqrt{\frac{z}{p}} \, e^{-\frac{\pi i}{12k}(h - h_0 p)} \, \eta\left( \frac{1}{k}\left(h + \frac{iz}{p} \right) \right), \label{poa1}
\end{align}
where
\begin{align}
w_{h,k} = \left( \frac{-h}{k} \right)
\exp\left( -\pi i \left[ \frac{k - 1}{4} + \frac{1}{12} \left( k - \frac{1}{k} \right)(h - h_0 p ) \right] \right). \label{poq1}
\end{align}

Similarly, applying~\eqref{poa} with $h^{\prime}=h_0$ and $w=-iz$ again, we find
\begin{align}
\eta\left( \frac{1}{k}\left(h_0 + \frac{i}{z} \right) \right)
&= w_{ph,k} \, e^{\frac{\pi i}{4}} \sqrt{z} \, e^{-\frac{\pi i}{12k}(ph - h_0)} \, \eta\left( \frac{1}{k}(ph + iz) \right), \label{poa2}
\end{align}
with
\begin{align}\label{poq2}
w_{ph,k} = \left( \frac{-ph}{k} \right)
\exp\left( -\pi i \left[ \frac{k - 1}{4} + \frac{1}{12} \left( k - \frac{1}{k} \right)(ph - h_0 ) \right] \right).
\end{align}

Combining equations~\eqref{poa1} and~\eqref{poa2} and taking into account \eqref{poq1} and \eqref{poq2}, we obtain
\begin{multline}
\eta^a\left( \frac{1}{k}\left(ph_0 +\frac{ip}{z} \right) \right) \cdot \eta^b\left( \frac{1}{k}\left(h_0 + \frac{i}{z} \right) \right)= \left( \frac{p}{k} \right)^b z^{\frac{a + b}{2}} p^{-\frac{a}{2}} e^{-\pi i \frac{(a + b)(k - 2)}{4}}\\
\times e^{ -\frac{\pi i}{12} \left\{ \left(k - \frac{1}{k} \right) h(a + bp) - k h_0(ap + b) \right\}} e^{-\frac{\pi z}{12pk}(a + pb)} \left( 1 + \sum_{n \geq 1} f_n q^n_{kp} \right),\label{on1}
\end{multline}
for some $f_n\in\mathbb{C}$, where $q_r:=q^{\frac{1}{r}}$.

\noindent{\bf Case 2: $p \mid k$.}\\
Let $k \in \mathbb{N}$ such that $p \mid k$, and let $H_0 \in \mathbb{N}$ satisfy $\gcd(H_0, k) = 1$. Set $h' = H_0$ and $w = z$. Choose $H$ such that
\[
H H_0 \equiv -1 \pmod{k}.
\]

Then, applying~\eqref{poa}, we find
\begin{align}\label{poa11}
\eta\left( \frac{1}{k/p}\left(h_0 +\frac{i}{z} \right) \right)
&= w_{H, k/p} \, e^{\frac{\pi i}{4}} \sqrt{z} \, e^{- \frac{\pi i}{12(k/p)}(H - H_0)} \, \eta\left( \frac{1}{k/p}(H + iz) \right),
\end{align}
with
\begin{align}\label{poq1a}
w_{H, k/p} = \left( \frac{-H}{k/p} \right)
\exp\left( -\pi i \left[ \frac{k/p - 1}{4} + \frac{1}{12} \left( \frac{k}{p} - \frac{p}{k} \right)(H - H_0 ) \right] \right).
\end{align}

Also,
\begin{align}\label{poa12}
\eta\left( \frac{1}{k}\left(H_0 + \frac{i}{z} \right) \right)
= w_{H,k} \, e^{\frac{\pi i}{4}} \sqrt{z} \, e^{- \frac{\pi i}{12k}(H - H_0)} \, \eta\left( \frac{1}{k}(H + iz) \right),
\end{align}
with
\begin{align}\label{poq12}
w_{H,k} = \left( \frac{-H}{k} \right)
\exp\left( -\pi i \left[ \frac{k - 1}{4} + \frac{1}{12} \left( k - \frac{1}{k} \right)(H - H_0) \right] \right).
\end{align}

Combining \eqref{poa11} and \eqref{poa12} and plugging in \eqref{poq1a} and \eqref{poq12}, we obtain
\begin{align}
&\eta^a\left( \frac{1}{k/p}\left(h_0 + \frac{i}{z} \right) \right) \cdot \eta^b\left( \frac{1}{k}\left(h_0 + \frac{i}{z} \right) \right) \notag \\
&= \left( \frac{-H}{k/p} \right)^a \left( \frac{-H}{k} \right)^b z^{\frac{a + b}{2}} \cdot
e^{ -\pi i \cdot \frac{(H - H_0)(a + b p)}{12k p} }
 e^{-\frac{\pi i}{4} \left[ \frac{k}{p}(a + b p) - (a + b) \right]}\nonumber\\
&\quad \cdot \exp\left( -\frac{\pi i}{12} \left[ \left( \frac{k}{p} - \frac{p}{k} \right)(aH - aH_0) + \left(k - \frac{1}{k} \right)(bH - bH_0) \right] \right) \notag\\
&\quad \cdot e^{-\frac{\pi z}{12k}(ap + b)}\left( 1 + \sum_{n \geq 1} l_n q^n_{k} \right)\label{on2}
\end{align}
for some $l_n\in\mathbb{C}$.



\vspace{.3cm}
\subsection*{Fourier Expansion of the \texorpdfstring{$\eta$}{eta}-Quotient under Hecke Operator}
To determine the growth of the numerator of \texorpdfstring{$g(q)$}{g(q)} on the right-hand side of \eqref{MM} towards the cusp \texorpdfstring{$0$}{zero}, we recall Proposition \eqref{prop1}
\begin{align}\label{UI}
&\left(\left(\eta|V_p\right)^{am+r} \eta^{mb-Mr}\right)\big|T(m)\left(\frac{-1}{\tau}\right)\nonumber \\
&=\left(\frac{p}{m}\right)^{am+r} m^{\kappa-1} \eta^{am+r}\left(\frac{-mp}{\tau}\right) \eta^{mb-Mr}\left(\frac{-m}{\tau}\right) \nonumber\\
&\quad +\frac{1}{m} \sum_{j=0}^{m-1} \eta^{am+r}\left( \frac{p(-1/\tau+j)}{m} \right) \eta^{mb-Mr}\left( \frac{-1/\tau + j}{m} \right)
\end{align}
where $\kappa$  is defined in \eqref{wt}.\\

\noindent \textbf{Case 1)} For $j \geq 1$: In our case, we let $\frac{\tau}{p}=iw = \frac{iz}{p}$, set $h_j' = pj$, and $k = m$. Then $h_j$ satisfies
\[
h_j \cdot pj \equiv -1 \pmod{m}, \quad 0 \leq h_j < m.
\]
For the term $\eta\left( \frac{p(-1/\tau + j)}{m} \right)$ in \eqref{poa}, we have:
\begin{align}\label{J111}
\eta\left( \frac{p(-1/\tau + j)}{m} \right)
= e^{\frac{\pi i}{4}} \sqrt{p^{-1}z} \cdot w_{h_j, m} \cdot e^{\frac{-\pi i}{12m}(h_j - pj)} \cdot \eta\left( \frac{h_j + izp^{-1}}{m} \right),
\end{align}
where
\begin{align}\label{w1111}
w_{h_j, m} = \left(\frac{-h_j}{m}\right) e^{-\pi i \left( \frac{1}{4}(m - 1) + \frac{1}{12}\left(m - \frac{1}{m}\right)(h_j - pj )\right)}.
\end{align}

Similarly,
\begin{align}\label{J222}
\eta\left( \frac{-1/\tau + j}{m} \right)
= e^{\frac{\pi i}{4}} \sqrt{z} \cdot w_{h_jp, m} \cdot e^{\frac{-\pi i}{12m}(h_jp - j)} \cdot \eta\left( \frac{ph_j + iz}{m} \right),
\end{align}
where
\begin{align}\label{w2222}
w_{h_jp, m} = \left(\frac{-h_jp}{m}\right) e^{-\pi i \left( \frac{1}{4}(m - 1) + \frac{1}{12}\left(m - \frac{1}{m}\right)(h_jp - j)\right)}.
\end{align}


Multiplying powers of \eqref{w1111} and \eqref{w2222}:
\begin{align}
w_{h_j, m}^{am + r} w_{h_jp, m}^{bm - Mr}
= \left(\frac{p}{m}\right)^{bm - Mr} e^{-\pi i \cdot \frac{m(a + b)(m - 1) + r(1 - M)}{4}} e^{-\frac{\pi i r}{12} \left(m - \frac{1}{m}\right)(h_j(1 - pM) + j(M - p))}.
\end{align}
Therefore, combining the above expression together with \eqref{J111} and \eqref{J222}, we obtain for $j \geq 1$,
\begin{align}
&\eta\left( \frac{p(-1/\tau + j)}{m} \right)^{am + r} \eta\left( \frac{-1/\tau + j}{m} \right)^{bm - Mr} \nonumber \\
&= \left(\frac{p}{m}\right)^{bm - Mr} p^{-\frac{am + r}{2}} z^{\frac{(a + b)m + r(1 - M)}{2}} \cdot e^{-\frac{\pi i r}{12} \left(m - \frac{1}{m}\right)(h_j(1 - pM) + j(M - p))} \nonumber \\
&\quad \times e^{-\frac{\pi i}{12m}(h_j - pj)(am + r)} e^{-\frac{\pi i}{12m}(h_jp - j)(bm - r)}  \times \eta\left( \frac{h_j + izp^{-1}}{m} \right)^{am + r} \eta\left( \frac{ph_j + iz}{m} \right)^{bm - Mr} \nonumber \\
&= \left(\frac{p}{m}\right)^{bm - Mr} i^{\frac{-[(a + b)m + r(1 - M)](m - 2)}{2}} p^{-\frac{am + r}{2}} z^{\frac{(a + b)m + r(1 - M)}{2}} \nonumber \\
& \qquad \qquad\times e^{\frac{\pi i r}{12m}(pM - 1)(m^2 - 1) h_j} e^{-\frac{\pi z}{12pm}(m(a + bp) + r(1 - Mp))} \left(1 + \sum_{n \geq 1} v_{n, m, j} q_{pm}^n \right):=A_j. \label{opj}
\end{align}

Recall $q_N = e^{\frac{2\pi i z}{N}}$. Summing over $j = 1$ to $m - 1$ and using von Sterneck’s formula:
\[
c_q(n) = \mu\left( \frac{q}{\gcd(q, n)} \right) \frac{\phi(q)}{\phi\left( \frac{q}{\gcd(q, n)} \right)}, \quad \text{where } c_q(n) = \sum_{\substack{1 \leq a \leq q \\ (a, q) = 1}} e^{2\pi i a n/q}.
\]
taking into account the fact $\gcd(r, m) = 1$ and $pM \not\equiv 1 \pmod{m}$, we get
\begin{align}
\sum_{j = 1}^{m - 1} A_j 
&= \left( \frac{p}{m} \right)^{bm - Mr} p^{-\frac{am + r}{2}} i^{-\kappa (m-2)}  z^{\kappa} \nonumber \\
&\quad \times e^{-\frac{\pi z}{12pm}(m(a + bp) + r(1 - Mp))} \left(-1 + \sum_{n \geq 1} \left( \sum_{j = 1}^{m - 1} v'_{n, m, j} \right) q_{pm}^n \right). \label{KL11}
\end{align}

\noindent\textbf{Case 2:} When $j=0$, we will apply the functional equation \eqref{fnal equation} for $\eta$-quotients to obtain (recall that $\tau=iz$)
\begin{align}
&\eta\left( \frac{-p/\tau}{m} \right)^{am+r} \eta\left( \frac{-1/\tau}{m} \right)^{bm - Mr}\nonumber\\
&= \left(\frac{zm}{p} \right)^{\frac{am + r}{2}} \left( -izm \right)^{\frac{bm - Mr}{2}} \eta\left( \frac{imz}{p} \right)^{am + r} \eta\left( imz \right)^{bm - Mr} \nonumber \\
&= p^{-\frac{am + r}{2}}  (zm)^{\kappa} e^{-\frac{\pi zm}{12p}(m(a + bp) + r(1 - pM))} \left(1 + \sum_{n \geq 1} c_n q_{pm}^{m^2 n} \right) := A_0. \label{opogeneral}
\end{align}

For the $V$-operator part, again applying the functional equation \eqref{fnal equation}:
\begin{align}
&\eta\left( \frac{-mp}{\tau} \right)^{am + r} \eta\left( \frac{-m}{\tau} \right)^{bm - Mr}\nonumber\\
&= p^{-\frac{am + r}{2}} \left(\frac{z}{p}\right)^{\kappa} e^{-\frac{\pi z}{12pm}(m(a + bp) + r(1 - pM))} \left(1 + \sum_{n \geq 1} s_n q_m^n \right). \label{r3general}
\end{align}

Combining \eqref{opogeneral}, \eqref{r3general}, and \eqref{KL11}, and keeping in mind that $M$ is odd, we obtain
\begin{align}
&\left(\left(\eta|V_p\right)^{am + r} \eta^{bm - Mr} \right)\big| T(m)\left(\frac{-1}{\tau}\right) \nonumber \\
&= \frac{1}{m} \left( \sum_{j = 1}^{m - 1} A_j + A_0 \right) + \frac{1}{m} \left( \frac{p}{m} \right)^{am + r} (-i)^{\frac{(a + b)m}{2}} p^{-\frac{am + r}{2}} z^{\frac{(a + b)m+r(1-pM)}{2}} \nonumber \\
&\quad \times e^{\frac{\pi iz}{12pm}(m(a + bp) + r(1 - p))} \left(1 + \sum_{n \geq 1} s_n q_m^n \right) \nonumber \\
&= \left( \frac{p}{m} \right)^{bm - Mr}  p^{-\frac{am + r}{2}} z^{\kappa} e^{-\frac{\pi  z}{12pm}(m(a + bp) + r(1 - pM))} \left( \sum_{n \geq 1} b_n q_{pm}^n \right).\label{Main11}
\end{align}
Since the numerator on the right hand side of \eqref{MM} belongs to $S_{\kappa}(\Gamma_0(p), \chi_p)$ by \eqref{modularform}, its Fourier expansion at $0$ can only have powers of $q_p$, so we can rewrite \eqref{Main11}, noting \eqref{on1}, as 
\begin{align}\label{main12}
g(q_1) =  z^{\kappa-(a+b)} c_p \cdot e^{-\frac{\pi z}{12pm}(r + p(24-Mr))} \sum_{n = 1}^\infty t_n q_p^n
\end{align}
for some $t_n\in\mathbb{C}$. Here $q_1=e^{2\pi i\frac{1}{k}(h'+\frac{i}{z})}\to 1$ as $q=e^{2\pi i \frac{1}{k}(h+iz)}\to 0$. Note that, for $m$ sufficiently large, we have 
\[
\frac{r}{24pm}+\frac{24-Mr}{24m}+ \frac{1}{p}>0,
\]
so \eqref{main12} decays as $q_1\to 0$.

Replacing $q$ by $q^{24/d}$ in \eqref{MM}, we obtain:
\begin{align}\label{main13}
\sum_{\substack{n \geq 0 \\ m \mid ((24/d) n + r(p - M)/d)}} c_{p, r, Mr}(n) q^{\frac{(24/d) n + r(p - M)/d}{m}} \equiv g\left(q^{24/d}\right) \pmod{m}.
\end{align}

Combining \eqref{main12}, \eqref{main13} and noting \eqref{MM}, we conclude:
\begin{equation}\label{eqn:gcuspidal}
g\left(q^{24/d}\right) \in S_{\kappa - \frac{(a + b)}{2}}\left( \Gamma_0\left( \left(\frac{24}{d} \right)^2 p \right), \chi_p \right).
\end{equation}
where we have used the fact $d = \gcd(p - M, \tfrac{24}{s})= \gcd(p M-1, \tfrac{24}{s})$.
Consequently, we will get
\[
\sum_{f \geq 0} c_{p, r, Mr}\left( \frac{dmf - r(p - M)}{24} \right) q^f \equiv \sum_{n \geq 0} u(n) q^n \pmod{m},
\]
where $g\left(q^{24/d}\right) = \sum_{n \geq 0} u(n) q^n$.

Now by Theorem \ref{thm:serre}, the set of primes $\ell$ such that
\[
\sum_{n = 0}^\infty u(n) q^n \big| T(\ell) \equiv 0 \pmod{m}
\]
has positive density, where $T(\ell)$ denotes the Hecke operator acting on
\[
S_{\kappa - \frac{(a + b)}{2}}\left( \Gamma_0\left( \left( \frac{24}{d} \right)^2 p \right), \chi_p \right).
\]
where $\kappa$ is defined in \eqref{wt}. Moreover, by the theory of Hecke operators, we have:
\[
\sum_{n = 0}^\infty u(n) q^n \big| T(\ell) = \sum_{n = 0}^\infty \left( u(\ell n) + \left( \frac{\ell}{p} \right) \ell^{\kappa-\frac{(a+b)}{2}-1} u\left( \frac{n}{\ell} \right) \right) q^n.
\]

Since $u(n)$ vanishes for non-integer $n$, we have $u\left(\frac{n}{\ell}\right) = 0$ when $(n, \ell) = 1$. Thus,
\[
u(\ell n) \equiv 0 \pmod{m} \quad \text{for } (n, \ell) = 1.
\]

Recalling that
\[
u(n) \equiv c_{p, r, Mr}\left( \frac{dmn - r(p - M)}{24} \right) \pmod{m},
\]
we obtain the final congruence:
\[
c_{p, r, Mr}\left( \frac{dmn\ell - r(p - M)}{24} \right) \equiv 0 \pmod{m}
\]
for each integer $n$ with $(n, \ell) = 1$.

\end{document}